\documentclass[twoside,a4paper,12pt] {article}
\usepackage{amsmath}
\usepackage{epsfig}\usepackage{psfrag}\usepackage{subfigure}
\usepackage{tabularx}\usepackage{rotating}
\newtheorem{theorem}{Theorem}

\newtheorem{lemma}{Lemma}

\makeatletter
\renewcommand{\@thesubfigure}{(\alph{subfigure})}
\renewcommand{\p@subfigure}{}
\renewcommand{\@makecaption}[2]{%
\vskip 10\p@   \setbox\@tempboxa\hbox{#1.\space#2}%
\ifdim \wd\@tempboxa >\hsize       #1.\space#2\par     \else       \hbox to\hsize{\hfil\box\@tempboxa\hfil}%
\fi}
\def\citefull{\def\astroncite##1##2{##1 ##2}\@internalcite}\def\@cite#1#2{#1\if@tempswa  #2\fi}

\def\@fnsymbol#1{\ifcase#1\or \mbox{${^{\star}}$}\or
   \dagger\or \ddagger\or
   \S \or \P \or \|\or \mbox{$^{\star\star}$}\or \dagger\dagger
   \or \ddagger\ddagger\or \S\S\or \P\P\or \|\|\else ***
   \fi\relax}
\newcommand{\lstar}{{\mbox{{\large $\star$}}}}
\makeatother 
\addtocounter{footnote}{0}
\begin{document}
\noindent {\Large \bf A Simple Proof of Inequalities of Integrals of Composite Functions} 
\thispagestyle{empty}
\vskip0.5cm
\noindent {\large Zhenglu Jiang$^{1\lstar}$\footnotetext{$^\lstar$E-mail: mcsjzl@mail.sysu.edu.cn (ZJ); 
mcsfxy@mail.sysu.edu.cn (XF); hongjiangtian@263.net (HT).}, 
Xiaoyong Fu$^{1\lstar}$ and Hongjiong Tian$^{2\lstar}$}
\newline{\footnotesize $^1$Department of Mathematics, Zhongshan University, 
Guangzhou 510275, China}
\newline{\footnotesize $^2$Department of Mathematics, Shanghai Teachers' University, 
Shanghai 200234, China}
\vskip0.5cm
\noindent{\footnotesize Submitted 2006 May 11} 
\baselineskip 16pt
\vskip0.6cm

\hfill\begin{minipage}{130mm}
\noindent{\bf  ABSTRACT} \newline
In this paper we give a simple proof of inequalities of integrals of functions   
which are the composition of nonnegative continous convex functions on a vector space ${\bf R}^m$ 
and vector-valued functions in a weakly compact subset  
of a Banach vector space generated by $m$ $L_\mu^p$-spaces for $1\leq p<+\infty.$ 
Also, the same inequalities hold if these vector-valued functions are in a  weakly* compact subset of   
  a Banach vector space generated by $m$ $L_\mu^\infty$-spaces instead. \newline 
\newline{\bf Key words:} convex functions, measure spaces, Lebesgue integrals, Banach spaces.
\end{minipage}
\vskip0.5cm

\section{Introduction}
\label{intro}
Both convexity of functions and characteristics of weakly compact sets  
are important in the study of extremum problems and integral estimates
in many areas of applied mathematics. The basic results about the 
continuity and differentiability of convex functions appear in the book 
of Rockafellar [\cite{rrt}] while  many characteristics 
of weakly compact sets in Banach spaces 
in the books of Benedetto [\cite{bjj}] and Yosida [\cite{yk}]. 
Of interest and important  is the study of integral estimates of a function 
which is the composition of a convex function on a vector space   
and a vector-valued function in a weakly compact subset in a Banach space  
since this kind of composite function often appears in many research fields such as 
compensated compactness methods (e.g., in [\cite{dl}] and [\cite{yl}]).  
Therefore it is very necessary to give the inequalities of integrals 
of the composite functions for our solving many problems in 
applied mathematics.  

Throughout this paper, ${\bf R}$ denotes the real number system, ${\bf R}^n$ is the usual vector 
space of real $n$-tuples $x=(x_1,x_2,\cdots,x_n),$ 
$\mu$  is a nonnegative Lebesgue measure of ${\bf R}^n,$ 
$L_\mu^p({\bf R}^n)$ represents a Banach space where each measurable function $u(x)$ 
has the following norm 
\begin{equation}
\parallel u(x) \parallel_p=\left(\int_{{\bf R}^n}|u(x)|^p d\mu\right)^{\frac{1}{p}}
\label{norm}
\end{equation} 
for any $p\in [1,+\infty),$  $(L_\mu^p({\bf R}^n))^m$ denotes a Banach vector space 
where each measurable vector-valued function has $m$ components in $L_\mu^p({\bf R}^n),$    
$L_\mu^\infty({\bf R}^n)$ represents a Banach space where each measurable function $u(x)$ 
has the following norm 
\begin{equation}
\parallel u(x) \parallel_\infty= 
\mathop{\hbox{ess}\sup}\limits_{x\in {\bf R}^n} |u(x)| 
\hbox{  } (\hbox{or say  } \parallel u(x) \parallel_\infty  
=\mathop{\inf\limits_ {E\subseteq {\bf R}^n}}
\limits_{\mu(E^C)=0}\max\limits_{x\in E} |u(x)| ) 
\label{norminf}
\end{equation} 
where $E^C$ represents the complement set of $E$ in ${\bf R}^n,$ 
 and  $(L_\mu^\infty({\bf R}^n))^m$ denotes a Banach vector space 
where each measurable vector-valued function has $m$ components in $L_\mu^\infty({\bf R}^n).$    
A sequence $\{u_i\}_{i=1}^{+\infty}$ is called to be weakly convergent to $u$ in 
$L_\mu^p({\bf R}^n)$  as $i\rightarrow +\infty$ for $1<p<+\infty,$ if  the following equality holds
\begin{equation}
 \mathop{\lim}\limits_{i\rightarrow+\infty}
\int_{{\bf R}^n}u_ivd\mu=\int_{{\bf R}^n}uvd\mu
\label{weqp}
\end{equation}
for all $v\in L_\mu^q({\bf R}^n)$ where $q=\frac{p}{p-1}$ and  $1<p<+\infty.$ 
A sequence $\{u_i\}_{i=1}^{+\infty}$ is called to be weakly convergent to $u$ in 
$L_\mu^1({\bf R}^n)$ as $i\rightarrow +\infty,$ if  the equality (\ref{weqp}) holds
for all $v\in L_\mu^\infty({\bf R}^n).$ 
A sequence $\{u_i\}_{i=1}^{+\infty}$ is called to be weakly* convergent to $u$ in 
$L_\mu^\infty({\bf R}^n)$ as $i\rightarrow +\infty,$ if  the equality (\ref{weqp}) holds
for all $v\in L_\mu^1({\bf R}^n).$ 
A sequence $\{u_i=(u_{1i},u_{2i},\cdots,u_{mi})\}_{i=1}^{+\infty}$ is called to 
be weakly convergent to $\hat{u}=(\hat{u}_{1},\hat{u}_{2},\cdots,\hat{u}_{m})$ in  
a Banach vector space $(L_\mu^p({\bf R}^n))^m$ as $i\rightarrow +\infty$ for $1\leq p<+\infty,$  
if  $\{u_{ji}\}_{i=1}^{+\infty}$ is weakly convergent to $\hat{u}_{j}$ 
in $L_\mu^p({\bf R}^n)$ as $i\rightarrow +\infty$ for all $j=1,2,\cdots,m$ and $1\leq p<+\infty.$
If  $\{u_{ji}\}_{i=1}^{+\infty}$ is weakly* convergent to $\hat{u}_{j}$ 
in $L_\mu^\infty({\bf R}^n)$ as $i\rightarrow +\infty$ for all $j=1,2,\cdots,m,$
a sequence $\{u_i=(u_{1i},u_{2i},\cdots,u_{mi})\}_{i=1}^{+\infty}$ is called to 
be weakly* convergent to $\hat{u}=(\hat{u}_{1},\hat{u}_{2},\cdots,\hat{u}_{m})$ in  
a Banach vector space $(L_\mu^\infty({\bf R}^n))^m$ as $i\rightarrow +\infty.$   

Let $f(x)$ be a function whose values are real or $\pm\infty$ and whose domain 
is a subset $S$ of ${\bf R}^m.$ $f(x)$ is called a convex function on $S$ if the set 
$\{(x,y)| x\in S, y\in {\bf R}, y\geq f(x)\}$ is convex as a subset of ${\bf R}^{m+1}.$ 
Then it is easily known that $f(x)$ is 
convex from $S$ to $(-\infty,+\infty]$ if and only if 
\begin{equation}
f(\lambda_1x_1+\lambda_2x_2+\cdots+\lambda_kx_k)\leq\lambda_1f(x_1)+\lambda_2f(x_2)
+\cdots+\lambda_kf(x_k)
\label{convex}
\end{equation}
whenever $S$ is a convex subset of ${\bf R}^m,$ $x_i  \in S$ ($i=1,2,\cdots$) 
$\lambda_1\geq 0,$ $\lambda_2\geq 0,$ $\cdots, $ $\lambda_k\geq 0,$ 
$\lambda_1+\lambda_2+\cdots+\lambda_k=1.$ 
This is called Jensen's inequality as $S={\bf R}^m. $ 

Inequalities of integrals of functions   
which are the composition of convex functions  
and vector-valued functions in a weakly* compact subset of $(L^\infty_\mu({\bf R}^n))^m$ are shown in [\cite{yl}] 
(See Theorem~\ref{th03} in Section~\ref{iw*cs}). 
In [\cite{dl}], if a special nonnegative convex function is considered 
and a weakly convergent sequence is constrained in $(L^1_\mu({\bf R}^n))^2,$ 
similar integrals of their composite functions are obtained 
(See Example 2 in Section~\ref{iwcs}). 
Therefore in this paper, by Fatou's lemma, we are going to give inequalities of integrals of functions   
which are the composition of nonnegative continous convex functions on a vector space ${\bf R}^m$ 
and vector-valued functions in a weakly compact subset of  
 a Banach vector space generated by $m$ $L_\mu^p$-spaces for $1\leq p<+\infty.$ 
Also, the same inequalities hold if these vector-valued functions are in a  weakly* compact subset of 
 a Banach vector space generated by $m$ $L_\mu^\infty$-spaces instead. 

The plan of this paper is as follows. In Section~\ref{iwcs}, we obtain 
inequalities of integrals of the composite functions $f(u)$ where 
$u$ is a limit of a weakly convergent sequence $\{u_i\}_{i=1}^{+\infty}$ in 
 $(L_\mu^p({\bf R}^n))^m$ for $1\leq p< +\infty$ and $f=f(x)$ is a nonnegative 
convex function from ${\bf R}^m$ to ${\bf R}.$ 
A similar result for a weakly* convergent sequence 
is shown in Section~\ref{iw*cs}.

\section{Inequalities for weakly convergent sequences}
\label{iwcs} 
The basic concepts have been introduced in the previous section.
In this section we show inequalities of integrals of functions   
which are the composition of nonnegative continous convex functions on a vector space $R^m$ 
and vector-valued functions in a weakly compact subset  of 
 a Banach vector space generated by $m$ $L_\mu^p$-spaces for $1\leq p<+\infty.$ 
That is the following
\begin{theorem}
Suppose that a sequence $\{u_i\}_{i=1}^{+\infty}$ 
 weakly converges in 
 $(L_\mu^p({\bf R}^n))^m$ to $u$ for $1\leq{p}<+\infty$ 
 as $i\rightarrow +\infty,$ where $m$ and $n$ are two positive integers.            
 If $f(x)$ is a nonnegative continous convex function from ${\bf R}^m$ to ${\bf R},$  
 then 
\begin{equation}
\mathop{\underline{\lim}}\limits_{i\rightarrow+\infty}
 \int_{\Omega}f(u_i)d\mu\geq\int_{\Omega}f(u)d\mu. \label{ineq}
\end{equation}
for all $\Omega\subseteq {\bf R}^n.$
\label{th01}
\end{theorem}

{\bf Remark 1}~~For $m=1,$ Theorem~\ref{th01} can be also written as follows. 
Suppose that a sequence $\{u_i\}_{i=1}^{+\infty}$ 
 weakly converges in 
 $L_\mu^p({\bf R}^n)$ to $u$ for $1\leq{p}<+\infty$ 
 as $i\rightarrow +\infty,$ where $n$ is a positive integers.            
 If $f(x)$ is a nonnegative continous convex function from ${\bf R}$ to ${\bf R},$  
 then the inequality (\ref{ineq}) holds. 

{\bf Remark 2}~~Let $S$ denote a convex subset of $ {\bf R}$ 
and $S^m$ the usual vector 
space of real $m$-tuples $x=(x_1,x_2,\cdots,x_m)$ 
where $x_i\in S$ ($i=1,2,\cdots,m$). Then $S^m$ is also convex. 
Suppose that $\{u_i\}_{i=1}^{+\infty}$ 
weakly converges in 
 $(L_\mu^p({\bf R}^n))^m$ to $u$ for $1\leq{p}<+\infty$ 
 as $i\rightarrow +\infty.$ 
Assume that all the values of $\{u_i\}_{i=1}^{+\infty}$ and of $u$ belong to $S^m$ and that 
$f(x)$ is a nonnegative continous convex function from $S^m$ to  ${\bf R}.$ 
Then the inequality (\ref{ineq}) also holds. 

{\bf Example 1}~~Assume that $\{u_i\}_{i=1}^{+\infty}$ 
is a nonnegative sequence which weakly  
converges in $L_\mu^1({\bf R}^n)$ to $u.$  Then 
\begin{equation}
\mathop{\underline{\lim}}\limits_{i\rightarrow+\infty}
 \int_{\Omega}u_id\mu\geq\int_{\Omega}ud\mu \label{ex01}
\end{equation}
for all $\Omega\subseteq {\bf R}^n.$ 
This is  obviously a special case of Remark~2 
when $m=1$ and $f(x)=x$ is defined in $[0,+\infty].$

{\bf Example 2}~~Step 2 of the theorem proof of DiPerna \& Lions 
in [\cite{dl}] shows the following result that 
$$\mathop{\underline{\lim}}\limits_{i\rightarrow+\infty}\int_{{\bf R}^n}
F(a_i,b_i)d\mu\geq\int_{{\bf R}^n}F(a,b)d\mu$$ if two positive sequences
 $\{a_i\}_{i=1}^{+\infty}$ and $\{b_i\}_{i=1}^{+\infty}$ 
  weakly converge in $L_\mu^1({\bf R}^n)$ to $a$ and $b$ respectively, where 
  $F(x,y)=(x-y)\log(\frac{x}{y})$ for $x>0$ and $y>0.$
Obviously, since $F(x,y)$ is a nonnegative continous convex 
 function defined in $(0,+\infty)\times(0,+\infty),$ 
DiPerna \& Lions' result is a special example of Remark~2
 for the case of $p=1$ and $m=2.$   
 
We are below going to give a 
  simple proof of Theorem~\ref{th01}.
In order to prove Theorem~\ref{th01}, we first recall the following result (See [\cite{bjj}]):

\begin{lemma}
Given a measure space $(X,{\cal A},\mu),$ $1\leq{p}<
+\infty,$ and $\{u_n,u:n=1,2,\cdots\}\subseteq L_\mu^p(X).$ Assume $u_n
\rightarrow{u}$ weakly. Then there is a subsequence $\{u_{n_k}:k=1,2,\cdots\}$ 
 whose arithmetic means $\frac{1}{m}\sum_{k=1}^mu_{n_k}$ converge in the 
 $L_\mu^p(X)$-topology to $u.$ 
\label{lemma01}
\end{lemma}

 Banach and Saks only proved Lemma~\ref{lemma01} for the $1<p<+\infty$ cases; the 
 result for $L_\mu^1(X)$ was showed by Szlenk in 1965. Using Lemma~\ref{lemma01},
   we can easily prove Theorem~\ref{th01}. 

{\bf Proof of Theorem~\ref{th01}}~~Put $\alpha_i=\int_{\Omega}f(u_i)d\mu$ ($i=1,2,\cdots $) and 
$\alpha=\mathop{\underline{\lim}}\limits_{i\rightarrow +\infty}\int_{\Omega}f(u)d\mu$ 
for all $\Omega\subseteq {\bf R}^n.$  
Then there exists a subsequence of $\{\alpha_i\}_{i=1}^{+\infty}$ such that 
this subsequence, denoted without loss of generality 
by $\{\alpha_i\}_{i=1}^{+\infty},$ converges to $\alpha$ as $i\rightarrow +\infty.$  

Since $u_i\rightarrow{u}$ weakly in $(L_\mu^p({\bf R}^n))^m$ for $1\leq p<+\infty,$ 
it is shown by Lemma~\ref{lemma01} that there exists a subsequence $\{u_{i_j}:j=1,2,\cdots\}$ 
such that $\frac{1}{k}\sum_{j=1}^ku_{i_j} \rightarrow u$ 
in $(L_\mu^p({\bf R}^n))^m$ for $1\leq p<+\infty$ as $k\rightarrow +\infty.$  
Thus there exists a subsequence of $\{\frac{1}{k}\sum_{j=1}^ku_{i_j}: k=1,2,\cdots\}$ 
such that this subsequence (also denoted without loss of generality by 
$\{\frac{1}{k}\sum_{j=1}^ku_{i_j}: k=1,2,\cdots\}$) satisfies that, as $k\rightarrow +\infty,$
\begin{equation}
\frac{1}{k}\sum_{j=1}^ku_{i_j} \rightarrow u  \hbox{  a.e. in  }{\bf R}^n.
\label{conv}
\end{equation} 
On the other hand, since $f(x)$ is a nonnegative continous convex 
function from ${\bf R}^m$ to ${\bf R},$ we have
\begin{equation}
f(\frac{1}{k}\sum_{j=1}^ku_{i_j})\leq \frac{1}{k}\sum_{j=1}^kf(u_{i_j}). 
\label{conf1}
\end{equation}
Hance, by (\ref{conf1}) and Fatou's lemma, we know that 
\begin{equation}
 \int_{\Omega}\mathop{\underline{\lim}}\limits_{k\rightarrow+\infty}
f(\frac{1}{k}\sum_{j=1}^ku_{i_j})d\mu
\leq\mathop{\underline{\lim}}\limits_{k\rightarrow+\infty}
\frac{1}{k}\sum_{j=1}^k\int_{\Omega}f(u_{i_j})d\mu. 
\label{conf2}
\end{equation}
It follows from (\ref{conv}) that 
\begin{equation}
 \int_{\Omega}f(u)d\mu
\leq\mathop{\underline{\lim}}\limits_{k\rightarrow+\infty}
\frac{1}{k}\sum_{j=1}^k\int_{\Omega}f(u_{i_j})d\mu, 
\label{conf3}
\end{equation}
or equivalently, 
\begin{equation}
 \int_{\Omega}f(u)d\mu
\leq\mathop{\underline{\lim}}\limits_{k\rightarrow+\infty}
\frac{1}{k}\sum_{j=1}^k\alpha_{i_j} 
\label{conf4}
\end{equation}
which gives (\ref{ineq}) since $\alpha_i\rightarrow\alpha$ 
as $i\rightarrow+\infty.$ This completes our proof. 

\section{Inequalities for weakly* convergent sequences}
\label{iw*cs}
In the previous section we have given 
inequalities of integrals of the composite functions for weakly convergent 
sequences in  $(L_\mu^p({\bf R}^n))^m$ for $1\leq p<+\infty.$ 
A similar result for weakly* convergent sequences 
in $(L_\mu^\infty({\bf R}^n))^m$ can be also below obtained  in this section.

Using the process of the proof of Theorem~\ref{th01}, we can prove the following theorem.

\begin{theorem}
Assume that a sequence $\{u_i\}_{i=1}^{+\infty}$ 
 weakly* converges in 
$(L_\mu^\infty({\bf R}^n))^m$ to $u$ as $i\rightarrow +\infty,$ 
where $m$ and $n$ are two positive integers. If $f(x)$ is a nonnegative 
continous convex function from ${\bf R}^m$ to ${\bf R},$  
then the inequality (\ref{ineq}) also holds 
for all $\Omega\subseteq {\bf R}^n.$
\label{th02}
\end{theorem}

{\bf Proof}~~Put $\Omega_R=\Omega\cap \{w: |w|<R,w\in {\bf R}^n\}.$ 
Then $\Omega_R$  is a bounded set 
in ${\bf R}^n$ for all the fixed positive real number $R.$ 
Since $u_i\rightarrow{u}$ weakly* in 
$(L_\mu^\infty({\bf R}^n))^m,$ $u_i\rightarrow{u}$ weakly* in 
$(L_\mu^\infty(\Omega_R)^m.$ Hance, 
by $L^\infty(\Omega_R)\subset L^1(\Omega_R),$ 
it can be easily known that 
$u_i\rightarrow{u}$ weakly in 
$(L_\mu^1(\Omega_R)^m.$ 
Then, using the process of the proof of Theorem~\ref{th01},
we can get 
\begin{equation}
\mathop{\underline{\lim}}\limits_{i\rightarrow+\infty}
 \int_{\Omega_R}f(u_i)d\mu\geq\int_{\Omega_R}f(u)d\mu. \label{ineq01}
\end{equation}
It follows from the nonnegativity of the convex function $f$ that 
\begin{equation}
\mathop{\underline{\lim}}\limits_{i\rightarrow+\infty}
 \int_{\Omega}f(u_i)d\mu\geq\int_{\Omega_R}f(u)d\mu. \label{ineq02}
\end{equation}
Finally, by Lebesgue dominated convergence theorem, 
as $R\rightarrow +\infty,$ (\ref{ineq02}) 
implies (\ref{ineq}).  Our proof is completed. 

Furthermore, in the case for continous convex functions, 
using Mazur's lemma [\cite{yk}],  we can deduce    
\begin{theorem}
Assume that a sequence $\{u_i\}_{i=1}^{+\infty}$ 
 weakly* converges in 
$(L_\mu^\infty({\bf R}^n))^m$ to $u$ as $i\rightarrow +\infty,$ 
where $m$ and $n$ are two positive integers. If $f(x)$ is a   
continous convex function from ${\bf R}^m$ to ${\bf R},$  
then the inequality (\ref{ineq}) also holds 
for all the bounded set $\Omega\subset {\bf R}^n.$
\label{th03}
\end{theorem} 

Obviously, the result in Theorem~\ref{th03} is weaker than that in Theorem~\ref{th02} 
since the nonnegativity of the convex function appears in Theorem~\ref{th02} but not in 
Theorem~\ref{th03}.  We can also deduce Theorem~\ref{th02} from Theorem~\ref{th03}. 
Theorem~\ref{th03} can be easily proved using the 
following lemma.
\begin{lemma}
 Assume $u_n\rightarrow{u}$ weakly in a normed linear space. 
Then there exists, for any $\epsilon>0,$ a convex combination 
$\sum_{k=1}^n\lambda_ku_{k}$ 
$(\lambda_k\geq 0, \sum_{k=1}^n\lambda_k=1)$ of $\{u_k:k=1,2,\cdots\}$ 
 such that $\parallel u-\sum_{k=1}^n\lambda_ku_{k} \parallel \leq \epsilon$ 
where $ \parallel v \parallel$ is a norm of $v$ in the space.
\label{lemma02}
\end{lemma}
This is called Mazur's lemma. Its proof can be found in the book of Yosida [\cite{yk}].  

In fact, Theorem~\ref{th03} is a part of the results given by Ying [\cite{yl}]. However, we 
still give its proof below. 

{\bf Proof of Theorem~\ref{th03}}~~Put $\alpha_i=\int_{\Omega}f(u_i)d\mu$ ($i=1,2,\cdots $) and 
$\alpha=\mathop{\underline{\lim}}\limits_{i\rightarrow +\infty}\int_{\Omega}f(u)d\mu$ 
for all the bounded set $\Omega.$  
Then there exists a subsequence of $\{\alpha_i\}_{i=1}^{+\infty}$ such that 
this subsequence, denoted without loss of generality 
by $\{\alpha_i\}_{i=1}^{+\infty},$ converges to $\alpha$ as $i\rightarrow +\infty.$  

Since $u_i\rightarrow{u}$ weakly* in 
$(L_\mu^\infty({\bf R}^n))^m,$  $u_i\rightarrow{u}$ weakly* in 
$(L_\mu^\infty(\Omega))^m.$ Hance, 
by $L^\infty(\Omega)\subset L^1(\Omega),$ 
it can be easily known that 
$u_i\rightarrow{u}$ weakly in 
$(L_\mu^1(\Omega)^m.$ 
It follows from Lemma~\ref{lemma02} that, 
for any natural number $j,$  there exists  a convex combination 
$\sum_{k=j}^{N(j)}\lambda_ku_{k}$ 
$(\lambda_k\geq 0, \sum_{k=j}^{N(j)}\lambda_k=1)$ of $\{u_k:k=j,j+1,\cdots\}$ 
 such that $\parallel u-\sum_{k=j}^{N(j)}\lambda_ku_{k} \parallel_\infty  
\leq \frac{1}{j}$ 
where $N(j)$ is a natural number which depends on $j$ and  
$\{u_k:k=j,j+1,\cdots\},$ $ \parallel v \parallel_\infty$ represents a norm of $v$ in $
(L_\mu^\infty(\Omega))^m.$ 
Put $v_{j}=\sum_{k=j}^{N(j)}\lambda_ku_{k}.$ Then, as $j\rightarrow +\infty,$ 
$v_j\rightarrow{u}$ in $(L_\mu^\infty(\Omega))^m.$ 
Since $f(x)$ is continous, for any given positive $\epsilon,$ 
there exists a natural number $N$ such as 
\begin{equation}
f(u)<f(v_j)+\frac{\epsilon}{\hbox{mes}(\Omega)}
\label{ineq04}
\end{equation}
for all $j>N,$ where $\hbox{mes}(\Omega)$ represents the measure of $\Omega.$ 
By the convexity of the function $f(x),$ integrating (\ref{ineq04}) gives
\begin{equation}
\int_{\Omega}f(u)d\mu\leq\sum_{k=j}^{N(j)}\lambda_k\int_{\Omega} f(u_{k})d\mu+\epsilon,
\label{ineq05}
\end{equation} 
or equivalently,
\begin{equation}
\int_{\Omega}f(u)d\mu\leq\sum_{k=j}^{N(j)}\lambda_k\alpha_k+\epsilon. 
\label{ineq06}
\end{equation} 
First let $j\rightarrow +\infty$ and then $\epsilon\rightarrow 0,$ 
(\ref{ineq06}) gives (\ref{ineq}) since $\alpha_k\rightarrow\alpha$ 
as $k\rightarrow +\infty.$  This completes the proof.

\section*{Acknowledgement}
This work was supported in part by NSFC 10271121. It was also sponsored by 
SRF for ROCS, SEM.


\end {document}